\documentclass[11pt,reqno]{amsart}
\usepackage{amssymb,latexsym}
\usepackage{graphicx}
\usepackage{mathtools}
\usepackage{color}
\usepackage[hyphens]{url}
\usepackage[T1]{fontenc}
\usepackage{hyperref}
\usepackage{amsmath}
\usepackage{mathdots}
\usepackage{breqn}
\usepackage[toc,page]{appendix}
\usepackage{url}    
\usepackage{breqn}
\usepackage{hyperref}
\usepackage{breakurl}
\newcommand{\bburl}[1]{\textcolor{blue}{\url{#1}}}

\calclayout

\makeatletter
\newcommand{\monthyear}[1]{%
  \def\@monthyear{\uppercase{#1}}}
\newcommand{\volnumber}[1]{%
  \def\@volnumber{\uppercase{#1}}}
\AtBeginDocument{%
\def\ps@plain{\ps@empty
  \def\@oddfoot{\@monthyear \hfil \thepage}%
  \def\@evenfoot{\thepage \hfil \@volnumber}}
\def\ps@firstpage{\ps@plain}
\def\ps@headings{\ps@empty
  \def\@evenhead{%
    \setTrue{runhead}%
    \def\thanks{\protect\thanks@warning}%
    \uppercase{\ }\hfil}%
  \def\@oddhead{%
    \setTrue{runhead}%
    \def\thanks{\protect\thanks@warning}%
    \hfill\uppercase{A Note on the Fibonacci Sequence and Schreier-type Sets}}%
  \let\@mkboth\markboth
  \def\@evenfoot{%
    \thepage \hfil \@volnumber}%
  \def\@oddfoot{%
    \@monthyear \hfil \thepage}%
  }%
\footskip=25pt
\pagestyle{headings}%
}
\makeatother

\theoremstyle{plain}
\numberwithin{equation}{section}

\newcommand{\seqnum}[1]{\href{https://oeis.org/#1}{\underline{#1}}}

\newcommand{\ignore}[1]{}

\newcommand\be{\begin{eqnarray}}
\newcommand\ee{\end{eqnarray}}
\newcommand\bea{\begin{eqnarray}}
\newcommand\eea{\end{eqnarray}}
\newcommand\ben{\begin{enumerate}}
\newcommand\een{\end{enumerate}}

\newtheorem{innercustomthm}{Theorem}
\newenvironment{customthm}[1]
  {\renewcommand\theinnercustomthm{#1}\innercustomthm}
  {\endinnercustomthm}

\newtheorem{innercustomcor}{Corollary}
\newenvironment{customcor}[1]
  {\renewcommand\theinnercustomcor{#1}\innercustomcor}
  {\endinnercustomcor}

\begin{document}

\monthyear{}
\volnumber{Volume, Number}
\setcounter{page}{1}
\title{A Note on the Fibonacci Sequence and Schreier-type Sets}

\author{H\`ung Vi\d{\^e}t Chu}

\address{Department of Mathematics, University of Illinois at Urbana-Champaign, Urbana, IL 61820} \email{hungchu2@illinois.edu}

\date{\today}

\begin{abstract}
A set $A$ of positive integers is said to be Schreier if either $A = \emptyset$ or $\min A\ge |A|$. 
We give a bijective map to prove the recurrence of the sequence $(|\mathcal{K}_{n, p, q}|)_{n=1}^\infty$ (for fixed $p\ge 1$ and $q\ge 2$), where 
$$\mathcal{K}_{n, p, q} \ = \ \{A\subset \{1, \ldots, n\}\,:\, \mbox{either }A = \emptyset \mbox{ or } (\max A-\max_2 A = p\mbox{ and }\min A\ge |A|\ge q)\}$$
and $\max_2 A$ is the second largest integer in $A$, given that $|A|\ge 2$. When $p = 1$ and $q=2$, we have that 
$(|\mathcal{K}_{n, 1, 2}|)_{n=1}^\infty$ is the Fibonacci sequence. As a corollary, we obtain a new combinatorial interpretation for the sequence $(F_n + n)_{n=1}^\infty$. 
\end{abstract}

\thanks{The author is thankful for the anonymous referee's suggestions that help improve the exposition of the note.}

\maketitle

A. Bird \cite{B} showed that for each $n\ge 1$, if we let
$$\mathcal{A}_n \ :=\ \{A\subset\{1, \ldots, n\}\,:\, n\in A\mbox{ and } \min A\ge |A|\},$$
then $|\mathcal{A}_n| = F_n$. The condition $\min A\ge |A|$ is called the \textit{Schreier condition}, and a set that satisfies the Schreier condition is called a \textit{Schreier set}. (The empty set satisfies the Schreier condition vacuously.) Schreier sets appeared in a paper of Schreier \cite{S} who used them to solve a problem in Banach space theory. The Schreier condition is also the central concept in a celebrated theorem by Odell \cite{Od}. Moreover, Schreier sets were independently discovered in combinatorics and appeared in Ramsey-type theorems for subsets of $\mathbb{N}$. Following the discovery by A. Bird, there has been research on various recurrences produced by counting Schreier-type sets (see \cite{BCF, C0, C1, C2, C3, M}). In this short note, we retrieve the Fibonacci sequence from a different counting problem than the one by A. Bird. In particular, for $n\ge 1$, define the set
$$\mathcal{K}_n\ :=\ \{A\subset [n]\,:\, \mbox{either }A = \emptyset \mbox{ or } (\max A-1\in A\mbox{ and }\min A\ge |A|)\},$$
where $[n] = \{1, \ldots, n\}$. While we fix the maximum element of sets in $\mathcal{A}_n$, we do not fix the maximum of sets in $\mathcal{K}_n$. Instead, we fix the distance between the largest and the second largest elements of sets in $\mathcal{K}_n$.

\begin{customthm}{A}\label{m1}
For $n\ge 1$, $|\mathcal{K}_n| = F_n$.
\end{customthm}

Let us briefly discuss the proof of Theorem \ref{m1}. It is easy to check that $|\mathcal{K}_1| = |\mathcal{K}_2| = 1$. We need only to show that $|\mathcal{K}_{n+1}| - |\mathcal{K}_n| = |\mathcal{K}_{n-1}|$ for all $n\ge 2$. Fix $n\ge 2$. By definition, $\mathcal{K}_n\subset \mathcal{K}_{n+1}$ and 
$$\mathcal{K}_{n+1}\backslash \mathcal{K}_n \ =\ \{A\subset [n+1]\,:\, n,n+1\in A\mbox{ and }\min A\ge |A|\}.$$  
We define a bijection $\pi: \mathcal{K}_{n-1}\rightarrow \mathcal{K}_{n+1}\backslash \mathcal{K}_n$: for $A\in \mathcal{K}_{n-1}$, 
$$\pi(A)\ :=\ \begin{cases}(A\backslash \{\max A\} + 1)\cup \{n, n+1\} &\mbox{ if } A\neq \emptyset,\\ \{n, n+1\}&\mbox{ if } A = \emptyset.\end{cases}$$
Interested readers may verify that $\pi$ is indeed a bijection or may look at the proof of the more general Theorem \ref{m2} below.

We have the following immediate corollary, which gives the sequence $(F_n + n)_{n=1}^\infty$ (see \seqnum{A002062}).
\begin{customcor}{B}
Let 
$$\mathcal{K}'_n\ :=\ \{A\subset [n]\,:\, \mbox{either } |A| \le 1 \mbox{ or } (\max A-1\in A\mbox{ and }\min A\ge |A|)\}.$$
Then $|K'_n| = F_n + n$ for all $n\ge 1$. 
\end{customcor}

\begin{proof}
Clearly, $|\mathcal{K}'_n| - |\mathcal{K}_n| = n$ for all $n\ge 1$. Using Theorem \ref{m1}, we are done. 
\end{proof}

We shall prove a more general result. Let $\max_2 A$ be the second largest number in $A$ if $|A|\ge 2$. For $n, p\ge 1$ and $q\ge 2$, define  
$$\mathcal{K}_{n, p, q} \ :=\ \{A\subset [n]\,:\, \mbox{either }A = \emptyset \mbox{ or } (\max A-\max_2 A = p\mbox{ and }\min A\ge |A|\ge q)\}.$$
Obviously, $\mathcal{K}_{n, 1, 2} = \mathcal{K}_n$.

\begin{customthm}{C}\label{m2}
Fix $n, p\ge 1$ and $q\ge 2$. We have
$$|\mathcal{K}_{n, p, q}| \ =\ \begin{cases}
1 &\mbox{ if } 1\le n\le p+2q-3,\\
|\mathcal{K}_{n-1, p, q}| + |\mathcal{K}_{n-2, p, q}| + \binom{n-p-q}{q-2} - 1 &\mbox{ if }n > p+2q-3.
\end{cases}$$
\end{customthm}

\begin{proof}
We prove Theorem \ref{m2} by recalling that $\mathcal{K}_{n-1, p, q}\subset \mathcal{K}_{n, p, q}$, then writing $$\mathcal{K}_{n, p, q}\backslash \mathcal{K}_{n-1, p, q} \ = \ \mathcal{S}\cup \mathcal{T}$$ for certain disjoint sets $\mathcal{S}$ and $\mathcal{T}$, and finally verifying that $|\mathcal{S}| = |\mathcal{K}_{n-2, p, q}|- 1$, while $|\mathcal{T}| = \binom{n-p-q}{q-2}$.

Fix $p\ge 1$ and $q\ge 2$. First, we check that for $1\le n\le p+2q-3$, $|\mathcal{K}_{n, p, q}| = 1$. Recall that
$$\mathcal{K}_{n, p, q} \ =\ \{A\subset [n]\,:\, \mbox{either }A = \emptyset \mbox{ or } (\max A-\max_2 A =  p\mbox{ and }\min A\ge |A|\ge q)\}.$$ Suppose $A$ is nonempty and $A\in \mathcal{K}_{n,p,q}$. Write $A = \{a_1, \ldots, a_k\}$, then $a_1\ge q$, $a_k\le p+2q-3$, and $a_{k-1}\le 2q-3$. Hence,
$$|\{a_1,\ldots, a_{k-1}\}|\ \le\ q-2$$
and so, $|A|\le q-1$, which contradicts the requirement that $|A|\ge q$. Therefore, for $1\le n\le p+2q-3$, $\mathcal{K}_{n,p,q} = \{\emptyset\}$. 

For $n\ge p+2q-2$, we shall show that $|\mathcal{K}_{n,p,q}|  = |\mathcal{K}_{n-1, p, q}| + |\mathcal{K}_{n-2, p, q}| + \binom{n-p-q}{q-2} - 1$. Let $\mathcal{S} = \{A\in \mathcal{K}_{n, p, q}\backslash \mathcal{K}_{n-1, p, q}: |A|\ge q+1\}$ and $\mathcal{T} = \{A\in \mathcal{K}_{n, p, q}\backslash \mathcal{K}_{n-1, p, q}: |A| = q\}$.
We define a bijection $\pi: \mathcal{K}_{n-2, p, q}\backslash \{\emptyset\}\rightarrow \mathcal{S}$: for a nonempty set $A\in \mathcal{K}_{n-2, p, q}$, 
$$\pi(A)\ :=\ (A\backslash \{\max A\} + 1)\cup \{n-p, n\}.$$

Firstly, $\pi$ is well-defined. Since $n\in \pi(A)$, $\pi(A)\notin \mathcal{K}_{n-1, p, q}$. That $\max A\le n-2$ implies that $\max_2 A\le n-2-p$, so $\pi(A)$ does not contain any number strictly between $n-p$ and $n$. Hence, $$\max \pi(A) - \max_2 \pi(A)\ =\ n - (n-p)\ =\ p.$$
Also, $|\pi(A)| = |A| + 1\ge q+1$ and 
$$\min \pi(A)\ =\ \min A + 1 \ \ge\ |A| + 1 = |\pi(A)|.$$
Therefore, $\pi(A)\in \mathcal{S}$.

Next, $\pi$ is one-to-one. Let $A_1, A_2\in \mathcal{K}_{n-2, p, q}\backslash \{\emptyset\}$ such that $\pi(A_1) = \pi(A_2)$. Note that $$\max (A_i\backslash \{\max A_i\} +1)\ \le\ (n-2-p)+1 \ =\ n-1-p,\mbox{ for }i = 1,2.$$
Hence, $\pi(A_1) = \pi(A_2)$ implies that $A_1\backslash \{\max A_1\} = A_2\backslash \{\max A_2\}$. So, $\max_2 A_1 = \max_2 A_2$, which, combined with $\max A_i - \max_2 A_i = p$ for $i=1, 2$, gives $A_1 = A_2$. We conclude that $\pi$ is one-to-one. 

Next, $\pi$ is onto. Take $A\in \mathcal{S}$. Then $n, n-p\in A$ and $|A|\ge q+1$. Let $B = A\backslash \{n-p, n\} - 1$ and $\ell = \max B$. Let $C = B\cup \{\ell + p\}$. We claim that $C\in \mathcal{K}_{n-2, p, q}$. Indeed, 
\begin{align*}\max C&\ =\ \max B + p\ \le\ n-p-1-1 + p \ =\ n-2,\\
\min C &\ =\ \min A - 1\ \ge\ |A| - 1\ =\ |B| + 1 \ =\ |C|,\mbox{ and }\\
|C| &\ =\ |B| + 1 \ =\ |A| - 1\ \ge\ (q+1)-1\ =\ q.\end{align*}
It is obvious from how we define $C$ that $\max C - \max_2 C = p$. Finally, $\pi(C) = A$ by construction. 

We have shown that $|\mathcal{S}| = |\mathcal{K}_{n-2, p, q}\backslash \{\emptyset\}| = |\mathcal{K}_{n-2, p, q}| - 1$. It remains to show that
$$|\mathcal{T}|\ =\ \binom{n-p-q}{q-2}.$$
A set $A$ is in $\mathcal{T}$ if and only if $\min A\ge |A| = q$, $\max A = n$, and $\max_2 A = n-p$. Hence, we can write 
a set $A$ in $\mathcal{T}$ as $A = D\cup \{n-p, n\}$, where $D\subset \{q, \ldots, n-p-1\}$ and $|D| = q-2$. Therefore, 
$|\mathcal{T}| = \binom{n-p-q}{q-2}$. This completes our proof as
\begin{align*}|\mathcal{K}_{n,p,q}|&\ =\ |\mathcal{K}_{n-1, p, q}| + |\mathcal{K}_{n,p,q}\backslash \mathcal{K}_{n-1, p, q}|\\
&\ =\ |\mathcal{K}_{n-1, p, q}| + |\mathcal{S}| + |\mathcal{T}|\\
&\ =\ |\mathcal{K}_{n-1, p, q}| + |\mathcal{K}_{n-2, p, q}| +  \binom{n-p-q}{q-2} - 1.
\end{align*}
\end{proof}

\ \\


\end{document}